\RequirePackage{ifpdf}
\ifpdf 
\documentclass[pdftex]{sigma}
\else
\documentclass{sigma}
\fi

\begin{document}

\allowdisplaybreaks

\renewcommand{\PaperNumber}{048}

\FirstPageHeading

\renewcommand{\thefootnote}{$\star$}

\ShortArticleName{Density of Eigenvalues of Random Normal Matrices}

\ArticleName{Density of Eigenvalues of Random Normal Matrices\\
with an Arbitrary Potential,\\ and of Generalized Normal Matrices\footnote{This paper
is a contribution to the Vadim Kuznetsov Memorial Issue
`Integrable Systems and Related Topics'. The full collection is
available at
\href{http://www.emis.de/journals/SIGMA/kuznetsov.html}{http://www.emis.de/journals/SIGMA/kuznetsov.html}}}

\Author{Pavel ETINGOF~$^*$ and Xiaoguang MA}

\AuthorNameForHeading{P. Etingof and X. Ma}

\Address{Department of Mathematics,
Massachusetts Institute of Technology, \\
77 Massachusetts Ave.,
Cambridge, MA 02139 USA}
\Email{\href{mailto:etingof@math.mit.edu}{etingof@math.mit.edu}, \href{mailto:xma@math.mit.edu}{xma@math.mit.edu}}
\URLaddressMarked{\url{http://www-math.mit.edu/~etingof/}}

\ArticleDates{Received December 05, 2006, in f\/inal form March
03, 2007; Published online March 14, 2007}

\Abstract{Following the works by Wiegmann--Zabrodin, Elbau--Felder,
Hedenmalm--Maka\-rov, and others, we consider the normal matrix
model with an arbitrary potential function, and explain how the problem of
f\/inding the support domain for the asymptotic eigenvalue density
of such matrices (when the size of the matrices goes to inf\/inity)
is related to the problem of Hele-Shaw f\/lows on curved surfaces,
considered by Entov and the f\/irst author in 1990-s. In the case
when the potential function is the sum of a rotationally
invariant function and the real part of a polynomial of the
complex coordinate, we use this relation and the conformal mapping method
developed by Entov and the f\/irst author to f\/ind the shape
of the support domain explicitly (up to f\/initely many
undetermined parameters, which are to be found from a f\/inite
system of equations). In the case when the rotationally invariant
function is $\beta |z|^2$, this is done by Wiegmann--Zabrodin and
Elbau--Felder.
We apply our results to the generalized normal matrix model,
which deals with random block matrices that give rise to
$*$-representations of the deformed preprojective algebra of the
af\/f\/ine quiver of type $\hat A_{m-1}$. We show that this model is
equivalent to the usual normal matrix model in the large $N$
limit. Thus the conformal mapping method
can be applied to f\/ind explicitly the support domain for the
generalized normal matrix model.}


\Classification{15A52} 

\renewcommand{\thefootnote}{\arabic{footnote}}
\setcounter{footnote}{0}

\section{Introduction}

The normal matrix model became
a focus of attention for many mathematical physicists
after the recent discovery (see e.g.~\cite{WZ1,KKMWZ,KMZ,MWZ}) of
its unexpected connections to the 2-dimensional dispersionless Toda
hierarchy and the Laplacian growth model (which is an exactly solvable model
describing free boundary f\/luid f\/lows in a Hele-Shaw cell or
porous medium). The original normal matrix model contained a
potential function whose Laplacian is a positive constant, but
later in \cite{WZ}, Wiegmann and Zabrodin considered
a more general model, where the potential function was
arbitrary. This is the model we will consider in this paper.

In the normal matrix model with an arbitrary potential function,
one considers the random normal matrices of some size $N$
with spectrum restricted to a compact domain $D$\footnote{The
compactness of $D$ is needed to guarantee convergence of the
arising integrals.} and probability measure
\begin{gather*}
 P_{N}(M)\text{d} M=Z_{N}^{-1}\exp(-N\textbf{tr} W(M))\text{d} M,
\end{gather*}
where $\text{d} M$ is the measure on the space of normal matrices
induced by the Euclidean metric on all complex matrices,
$W$ is a potential function (a real function on $\Bbb C$ with
some regularity properties, e.g. continuous),
and $Z_{N}$ is a normalizing factor.

In the original works on the normal matrix model,
the potential was
\begin{gather*}
W(z)=\beta |z|^2-P(z)-\overline{P(z)},
\end{gather*}
where $P$ is a complex polynomial of some degree $d$, and $\beta$
a positive real number. For this type of potential,
it was shown in the works \cite{WZ1,KKMWZ,KMZ,MWZ} (and then
proved rigorously in \cite{EF}) that under some conditions on the
potential, the asymptotic density of
eigenvalues is uniform with support in the
interior domain of a closed smooth curve. This curve is a
solution of an inverse moment problem, appearing in the theory of
Hele-Shaw f\/lows with a free boundary. Thus, applying the conformal
mapping method
(see \cite{VE} and references therein), one discovers that the conformal map
of the unit disk onto the outside of this curve which maps
$0$ to $\infty$ is a Laurent polynomial of degree $d$. This allows one to f\/ind
the curve explicitly up to f\/initely many parameters, which can be
found from a f\/inite system of algebraic equations.

In \cite{WZ}, Wiegmann and Zabrodin generalized this analysis
to an arbitrary potential function. They showed that
the density of eigenvalues is the Laplacian of the potential
function, and the eigenvalues are concentrated in the domain
which can be determined from an appropriate inverse moment
problem. This was proved rigorously in the paper \cite{HM}, which
extends the Elbau--Felder work to the case of an arbitrary
potential.

One of the goals of the present paper is
to use the generalized conformal mapping method,
developed in \cite{EE} by Entov and the f\/irst author
for studying Hele-Shaw f\/lows with moving boundary for curved
surfaces,
to calculate the boundary of the region of eigenvalues
explicitly in the case when
\begin{gather}\label{pote}
W(z)=\Phi(|z|^2)-P(z)-\overline{P(z)},
\end{gather}
where $\Phi$ is a function of one variable.
In this case, the conformal map of the disk onto the outside of
the curve is no longer algebraic, but one can still give an
explicit answer in terms of a contour integral.

Another goal is to extend the above results to the case
of {\it generalized normal matrix model}.
In this model, we consider block complex matrices of a certain kind
with commutatation relations similar to the def\/inition of a
normal matrix; they give rise to $*$-representations of the
deformed preprojective algebra of the
af\/f\/ine quiver of type $\hat A_{m-1}$. We prove that the problem of
computing the asymptotic eigenvalue distribution for this model, as the size of
the matrices goes to inf\/inity, is equivalent to the same problem
for the usual normal matrix model. This allows one to f\/ind
the boundary of the eigenvalue region explicitly if the potential
is given by (\ref{pote}).

The structure of this paper is as follows.
In Section \ref{NMM}, we state some basic facts
about the normal matrix model.
In Section \ref{BNMM}, we def\/ine the generalized normal
matrix model, and write down the probability measure in this
model. In Section \ref{EM}, we recall some facts about the
equilibrium measure and explain that the asymptotic eigenvalue
distribution tends to the equilibrium measure in the normal
matrix model and the generalized normal matrix model. In Section
\ref{RBD}, we use the singular point method from \cite{EE,VE} to
reconstruct the boundary of the support domain of the equilibrium
measure.

\section{The normal matrix model with an
arbitrary potential function}\label{NMM}

Let $D$ be a domain in the complex plane $ \mathbb{C}$.
We consider the set
\begin{gather*}
\mathcal{N}(D)=\{M\in \text{Mat}_{N}( \mathbb{C})\vert[M,M^\dagger]=0,
{\rm spectrum}(M)\subset D\}
\end{gather*}
of normal matrices with spectrum in the
domain $D$. Let $\text{d} M$ be the measure on $\mathcal{N}(D)$
induced by the Euclidean metric on $\text{Mat}_{N}( \mathbb{C})$.
It is well known (see e.g.~\cite{O,CZ}), that in terms of the
eigenvalues this measure on $\mathcal{N}(\Bbb C)$ is given by the formula
\begin{gather*}
\text{d} M=\prod_{1\leq i<j\leq N}|z_i-z_j|^2\prod_{i=1}^{N}\text{d}^2
z_i \text{d} U,
\end{gather*}
where $M=U\textbf{diag}(z_1,\ldots,z_N)U^{\dagger}$, $U\in U(N)$,
and $\text{d} U$ denotes the normalized $U(N)$-invariant measure on
the f\/lag manifold $U(N)/U(1)^N$.

Now let $W: \Bbb C\to \Bbb R$ be a continuous function.
If $M$ is a normal matrix, then we can def\/ine
$W(M)$ to be $\textbf{diag}(W(z_1),\ldots,W(z_N))$ in an orthonormal basis
in which $M=\textbf{diag}(z_1,\ldots,z_N)$.
It follows from the above that
the probability measure on $\mathcal{N}(D)$
with potential function $W$ is given by
\begin{gather}\label{NM-PF}
 P_N(M)\text{d} M=Z_N^{-1}e^{-N\sum_i W(z_i)}
\prod_{1\leq i<j\leq N}|z_i-z_j|^2\prod_{i=1}^{N}\text{d}^2 z_i \text{d} U,
\end{gather}
where $Z_N=\int_{D^N}e^{-N\sum_i
W(z_i)}\prod\limits_{1\leq i<j\leq N}|z_i-z_j|^2\prod\limits_{i=1}^{N}\text{d}^2 z_i$.
Here we assume that the integral is convergent (this is the case,
for instance, if $D$ is compact).

\section{The generalized normal matrix model}\label{BNMM}

\subsection{Generalized normal matrices}
Let us consider the following
generalization of normal matrices. Let $m\ge 1$ be an integer.
For a f\/ixed collection $\lambda=(\lambda_1,\ldots,\lambda_m)$
of real numbers such that $\sum_i \lambda_i=0$, and a domain $D$, we
def\/ine $\mathcal{N}_m(\lambda,D)$ to be the subset of $A\in \text{Mat}_{mN}( \mathbb{C})$
satisfying the following conditions:
\\ \indent for any $A\in \mathcal{N}_m(\lambda,D)$,
\begin{itemize}
\itemsep=0pt
 \item If $A_{ij}$,
$i,j=1,\ldots,m$ are $N\times N$ blocks of $A$, then $A_{ij}=0$ unless $j-i=1$ mod $m$;
 \item The spectrum of $A_{12}A_{23}\cdots A_{m1}$ is contained in $D$;
 \item $[A,A^\dagger]=\left(\begin{array}{lll}
\lambda_1I_N &  & \\
 & \ddots &  \\
 &  & \lambda_mI_N\end{array}\right)
$, where $I_N$ is the identity matrix of size $N$.
\end{itemize}

Note that $\mathcal{N}_1(0,D)=\mathcal{N}(D)$, thus elements
of $\mathcal{N}_m(\lambda,D)$ are a generalization of normal matrices.
We will thus call them generalized normal matrices.

\begin{remark} Generalized normal matrices are related
in the following way to quiver representations.
Let $Q$ be the cyclic quiver of type $\hat A_{m-1}$,
and $\bar Q$ its double. Let $\Pi_Q(\lambda)$ be the
deformed preprojective algebra of $Q$ with parameters $\lambda$
(see \cite{CBH}). By def\/inition, this algebra is the quotient of
the path algebra of $\bar Q$ by the relation
$\sum\limits_{a\in Q} [a,a^*]=\sum \lambda_i e_i$, where $e_i$ are the
vertex idempotents. The algebra $\Pi_Q$ has a $*$-structure,
preserving $e_i$ and sending $a$ to $a^*$ and $a^*$ to $a$.
It is easy to see that $\mathcal{N}_m(\lambda,D)$ is the set of all
matrix $*$-representations of $\Pi_Q$ of dimension $N\delta$
(where $\delta=(1,1,\ldots,1)$ is the basic imaginary root)
such that the spectrum of the monodromy operator
$a_1\cdots a_m$ is in $D$.
\end{remark}

Denote $A_{i,i+1}$ by $A_i$.
The group $U(N)^m=\underbrace{U(N)\times U(N)\times\cdots\times
U(N)}_{m \text{\ times\ }}$ acts naturally on $\mathcal{N}_m(\lambda,\Bbb
C)$ by the formula
\begin{gather*}
 (S_1,\ldots,S_m)A=\left(\begin{array}{ccccc} &
S_1A_1S_2^\dagger &  &  &  \\ &  & S_2A_2S_3^\dagger &  &  \\ &  &  & \ddots  &  \\ &  &  &  & S_{m-1}A_{m-1}S_m^\dagger \\S_mA_mS_1^\dagger &  &  &  & \end{array}\right).
\end{gather*}

We have the following lemma, which is a generalization
of the fact that a normal matrix diagonalizes in an orthonormal basis:
\begin{lemma}\label{diagonal}
 For any element $A\in \mathcal{N}_m(\lambda, D)$, we can find an element $ (S_1,\ldots,S_m)\in U(N)^m$ such that
 \begin{gather*}
  (S_1,\ldots,S_m)A=\left(\begin{array}{ccccc} & D_1 &  &  &  \\ &  & D_2 &  &  \\ &  &  & \ddots  &  \\ &  &  &  & D_{m-1} \\D_{m}  &  &  &  & \end{array}\right),
  \end{gather*}
   where $D_i$ are diagonal matrices.
\end{lemma}

\begin{proof}
 From the def\/inition, we have $A_iA_i^\dagger-A_{i-1}^\dagger
A_{i-1}=\lambda_i I_{N}$, where the index is considered modulo $m$,
and $\sum\limits_{i=1}^m \lambda_i=0.$

Now consider a collection of $N$-dimensional
unitary spaces $\{V_i\}_{i=1}^m$, and let us regard
$A_i$ as a~linear map $A_i: V_{i+1}\to V_i.$
So $A_iA_i^\dagger$ is a hermitian endomorphism of $V_i$.

Now suppose that $v$ is an eigenvector of
$A_iA_i^\dagger$ with eigenvalue $\nu$. We claim
that then $A_{i-1}v$ (if it is nonzero) is an eigenvector of
$A_{i-1}A_{i-1}^\dagger$ with eigenvalue $\nu-\lambda_i$.
Indeed,
\begin{gather*}
A_{i-1}A_{i-1}^\dagger A_{i-1}v=A_{i-1}(A_iA_i^\dagger-\lambda_i)v
=(\nu-\lambda_i)A_{i-1}v.
\end{gather*}
Thus, denoting by $V_i(\nu)$ the eigenspace of $A_iA_i^\dagger$
with eigenvalue $\nu$, we f\/ind that $A_{i-1}: V_i(\nu)\to
V_{i-1}(\nu-\lambda_i)$. Since $V_i=\oplus_{\nu\in \Bbb R}
V_i(\nu)$ (as $A_iA_i^\dagger$ is hermitian),
it suf\/f\/ices to prove the lemma in the case when $A_iA_i^\dagger$
is a scalar in $V_i$, in which case the statement is easy.
\end{proof}

\subsection{The Euclidean measure on generalized normal matrices}

First, let us consider the $N=1$ case.
Pick real numbers $\alpha_i$ such that
$\lambda_i=\alpha_i-\alpha_{i-1}$, and let $Q(x)=\prod\limits_{i=1}^m
(x+\alpha_i)$. A matrix
\begin{gather*}
A=\left(\begin{array}{ccccc} & r_1e^{i\theta_1} &  &  &  \\ &  & r_2e^{i\theta_2} &  &  \\ &  &  & \ddots  &  \\ &  &  &  & r_{m-1}e^{i\theta_{m-1}} \\r_me^{i\theta_m} &  &  &  & \end{array}\right),
\end{gather*}
(where $r_j\ge 0$, $\theta_j\in [0,2\pi)$)
is in $\mathcal{N}_m(\lambda,\Bbb C)$ if and only if
\begin{gather*}
r_{1}^{2}-r_{m}^{2} = \lambda_{1} = \alpha_{1}-\alpha_{m},\\
r_{2}^{2}-r_{1}^{2}=\lambda_{2}=\alpha_{2}-\alpha_{1},\\
\qquad \vdots\\
r_{m}^{2}-r_{m-1}^{2}=\lambda_{m}=\alpha_{m}-\alpha_{m-1}.
\end{gather*}

Thus to each $A\in \mathcal{N}_m(\lambda,\Bbb C)$ we can attach a real number
$x=r_{i}^{2}-\alpha_{i}$, which is independent of $i$, and
a complex number $z=\prod\limits_{j=1}^m r_j e^{i\theta_j}$.
It is easy to see that the point $(z,x)$ belongs to the surface
\begin{gather*} 
\Sigma=\Big\{(z,x)\in  \mathbb{C}\times \mathbb{R} \,|\, x+\alpha_i\ge 0\ \forall\, i,\
z\bar{z}=\prod_{i=1}^{m}(x+\alpha_{i})\Big\}.
\end{gather*}

Moreover, it is clear that any point of $\Sigma$ corresponds to
some $A$, and two matrices $A$, $A'$ giving rise to the same point
$(z,x)$ are conjugate. This implies that we have a
bijection between the equivalence classes in $\mathcal{N}_m(\lambda,\Bbb C)$
under the action of $U(1)^m$ and points of $\Sigma$.
Writing $z=re^{i\theta}$, we see that $x$, $\theta$ are coordinates
on $\Sigma$, so we may write the Euclidean measure on
$\mathcal{N}_m(\lambda,\Bbb C)$ using the coordinates $x$, $\theta$.

\begin{theorem}\label{eucmes1}
The Euclidean measure on $\mathcal{N}_m(\lambda,\Bbb C)$ for $N=1$ is:
 \begin{gather*}
\mu=\frac{1}{2}\sqrt{Q'(x)}\text{\rm d} x\text{\rm d}\theta \text{\rm d} U,
\end{gather*}
where $\text{\rm d} U$ is the Haar measure on $U(1)^m/U(1)$.
\end{theorem}

\begin{proof}
We have $r_{i}=\sqrt{x+\alpha_{i}}$.
Thus the parametrized real curve
$\{(r_{1}(x),\ldots,r_{m}(x))|x\in\mathbb{R}\}$ has length element
\begin{gather*}
\text{d} s=\sqrt{\sum_{i}r_{i}^{'}(x)^{2}}\text{d} x =
\frac{1}{2}\sqrt{\sum_i\frac{1}{x+\alpha_i}}\text{d} x.
\end{gather*}
This implies that the Euclidean measure on $\mathcal{N}_m(\lambda,\Bbb C)$ is
\begin{gather*}
\text{d}\mu
=\frac{1}{2}\sqrt{\sum_i\frac{1}{x+\alpha_i}}\text{d} x\prod_{j}r_{j}\text{d}\theta_{j}
=\frac{1}{2}\sqrt{\sum_i\prod_{j\neq
i}(x+\alpha_j)}\text{d} x\text{d}\theta_{1}\cdots \text{d}\theta_{m}\\
\phantom{\text{d}\mu}{}=\frac{1}{2}\sqrt{Q'(x)}\text{d} x\text{d}\theta_{1}\cdots \text{d}\theta_{m}
=\frac{1}{2}\sqrt{Q'(x)}\text{d} x\text{d}\theta \text{d} U,
\end{gather*}
as desired.
\end{proof}

Let us now consider the case of general $N$.
From Lemma \ref{diagonal}, we know that under the action of
$U(N)^m$, the equivalence class of $A\in \mathcal{N}_m(\lambda,\Bbb C)$
can be represented by $m$ diagonal matrices{\samepage
 \begin{gather}\label{equ-class}
  M=\left(\begin{array}{ccccc} & D_1 &  &  &  \\ &  & D_2 &  &  \\ &  &  & \ddots &  \\ &  &  &  & D_{m-1} \\D_m &  &  &  & \end{array}
\right),
 \end{gather}
where
 $D_i=\textbf{diag}(z_{1}^{i},\ldots,z_{N}^{i})$.}

From the def\/inition
\begin{gather*}
D_{1}D_{1}^{\dagger}-D_{m}^{\dagger}D_{m}=\lambda_{1}I_{N},\\
D_{2}D_{2}^{\dagger}-D_{1}^{\dagger}D_{1}=\lambda_{2}I_{N},\\
\qquad \vdots\\
D_{m}D_{m}^{\dagger}-D_{m-1}^{\dagger}D_{m-1}=\lambda_{m}I_{N}.
\end{gather*}

So we have
$z^{i}_{j}\overline{z^{i}_{j}}-z^{i+1}_{j}\overline{z^{i+1}_{j}}=\alpha_{i}-\alpha_{i+1}$.
Let $x_{j}=z^{i}_{j}\overline{z^{i}_{j}}-\alpha_{i}$ and $z_{j}=
\prod_{i}z^{i}_{j}$, then we have
\begin{gather*}
|z_{j}|^{2}=\prod_{i} (x_{j}+\alpha_{i}), \qquad j=1,\ldots,N.
\end{gather*}
Thus $((z_1,x_1),\ldots,(z_N,x_N))$ is a point on $\Sigma^N/S_N$.

Similarly to $N=1$ case, it is easy to show that
this gives rise to a bijection between
conjugacy classes of elements of $\mathcal{N}_m(\lambda,\Bbb C)$
and points of $\Sigma^N/S_N$. Using this fact
and combining the method of computation for usual normal matrices
with the $N=1$ case, one gets the following result.

\begin{theorem}\label{eucmes}
The Euclidean measure on $\mathcal{N}_m(\lambda,\Bbb C)$ has the form
\begin{gather*}
\text{\rm d} M=\frac{1}{2^{N}}\prod_{i}\sqrt{Q'(x_{i})}
\prod_{i<j} |z_{i}-z_{j}|^{2}\text{\rm d} x_{1}\cdots \text{\rm d} x_{N}
\text{\rm d} \theta_{1}\cdots\text{\rm d} \theta_{N}\text{\rm d} U,
\end{gather*}
where $\text{\rm d} U$ is the normalized invariant measure on
$U(N)^m/U(1)^N$.
\end{theorem}
\begin{proof}
At f\/irst, consider the subset $\mathcal{N}_m^{\rm diag}(\lambda,\Bbb C)$ of
$\mathcal{N}_m(\lambda,\Bbb C)$
consisting of the elements $M$ of the form  (\ref{equ-class}).
Then by Theorem \ref{eucmes1},
the measure on $\mathcal{N}_m^{\rm diag}(\lambda,\Bbb C)$ induced by the Euclidean metric
 is the product measure:
\begin{gather}\label{dia}
\mu_{\rm diag}=\frac{1}{2^{N}}\prod_{i}\sqrt{Q'(x_{i})}\text{d} x_{1}\cdots \text{d} x_{N}
\text{d} \theta_{1}\cdots\text{d} \theta_{N}\text{d} U_{\rm diag},
\end{gather}
where $\text{d} U_{\rm diag}$ is the Haar measure on
$U(1)^{Nm}/U(1)^N$.

Now consider the contribution of the of\/f-diagonal part.
Consider the elements
\begin{gather*}
\{v_{i,j}=E_{i,j}-E_{j,i},w_{i,j}=\sqrt{-1}(E_{i,j}+E_{j,i})\, |\, 0\leqslant i<j\leqslant N\}
\end{gather*}
of the Lie algebra of $U(N)$.

Let $V_{i,j,k}$, $W_{i,j,k}$ be the derivatives of
$(\exp(tv_{i,j}))_kM$ and $(\exp(tw_{i,j}))_kM$ at $t=0$, where
$a_k:=(1,\dots,1,a,1,\dots,1)\in U(N)^m$, with $a\in U(N)$ in the $k$-th place.
Then by formula (\ref{dia}), we have
\begin{gather*}
\text{d} M=\phi\cdot \frac{1}{2^{N}}\prod_{i}\sqrt{Q'(x_{i})}\text{d} x_{1}\cdots \text{d} x_{N}
\text{d} \theta_{1}\cdots\text{d} \theta_{N}\text{d} U,
\end{gather*}
where
\begin{gather}\label{phi}
\phi=|\wedge_{i<j,k}(V_{i,j,k}\wedge W_{i,j,k})|.
\end{gather}

To calculate $\phi$, let us denote by $B_{i,j,k}$, $i\ne j$,
the derivative of $(\exp(tE_{i,j}))_kM$ (note that since
$E_{i,j}$ lies only in the complexif\/ied Lie algebra
of $U(N)^m$, we have $(\exp(tE_{i,j}))_kM\notin \mathcal{N}_m(\lambda,\Bbb C)$,
but this is not important for our considerations). Then
equation (\ref{phi}) takes the form
\begin{gather*}
\phi=|\wedge_{i\ne j,k}B_{i,j,k}|.
\end{gather*}

Now $\phi$ can be easily calculated. To do so, we note that for a
given $i$, $j$, the transformation $(\exp(tE_{i,j}))_k$ changes
only the entries $a_{i,j}^p$ of $M$. On these entries, it acts
by
\begin{gather*}
a_{i,j}^p\to a_{i,j}^p+t(z_j^p\delta_{p,k}-z_i^p\delta_{p,k-1}).
\end{gather*}
This means that for each $i$, $j$,
$|\wedge_k B_{i,j,k}|=|J_{i,j}|$, where
\begin{gather*}
J_{i,j}:=\det\left(\begin{array}{ccccc}z^{1}_{j} &- z_{i}^1 & 0 & \cdots & 0 \\0 & z^{2}_{j} &- z_{i}^2 & \ddots & \vdots \\\vdots & 0 &z^{3}_{j} & \ddots & 0 \\0 & \ddots & 0 & \ddots &- z_{i}^m \\-z_{i}^m & 0 & \cdots & 0 & z^{m}_{j}\end{array}\right)=\prod_{s=1}^{m}z_{j}^{s}-\prod_{s=1}^{m}z_{i}^{s}=z_{j}-z_{i}.
\end{gather*}
This implies that
\begin{gather*}
\phi=\prod_{i\ne j}|J_{i,j}|=\prod_{i<j} |z_{i}-z_{j}|^{2},
\end{gather*}
as desired.
\end{proof}

\subsection{The probability measure with potential function\\
on generalized normal matrices}\label{Prob-Mea-PF}

Let $W:\Bbb C\to \Bbb R$ be a potential function.
The probability measure on $\mathcal{N}_m(\lambda,D)$ corresponding to
this function is def\/ined similarly to the case of usual normal matrices:
\begin{gather*}
 P_{N}(M)\text{d} M=Z_{N}^{-1}\exp(-N\textbf{tr} W(M_1\cdots M_m))\text{d} M,
\end{gather*}
$M\in \mathcal{N}_m(\lambda,D)$, where $M_i$ are the blocks of $M$.
Thus in terms of eigenvalues
\begin{gather*}
 P_{N}(M)\text{d} M=
\frac{1}{2^{N}Z_N}\exp\left\{-N\sum_j W(z_j)\right\}\\
\phantom{P_{N}(M)\text{d} M=}{}\times
\prod_{i}\sqrt{Q'(x_{i})}
\prod_{i<j} |z_{i}-z_{j}|^{2}\text{d} x_{1}\cdots \text{d} x_{N}
\text{d} \theta_{1}\cdots\text{d} \theta_{N}\text{d} U.
\end{gather*}

\begin{example}\label{mmdag} Let us calculate the potential
function corresponding to the quadratic potential ${\rm Tr}(MM^\dagger)$.
We have
\begin{gather*}
{\rm Tr}(MM^\dagger)=\sum_{i,j}|z_j^i|^2=\sum_{i,j}(x_j+\alpha_i)=
N\sum_i\alpha_i+m\sum_j x_j.
\end{gather*}
Thus if we choose $\alpha_i$ so that $\sum_i \alpha_i=0$
(this can be done in a unique way), then
\begin{gather*}
{\rm Tr}(MM^\dagger)=m\sum_j x_j,
\end{gather*}
so the corresponding potential function is
$W(z)=mQ^{-1}(|z|^2)$ (the function $Q$ is invertible on the
interval $[-\alpha,\infty)$, where $\alpha=\min \alpha_i$).
\end{example}

\section{Equilibrium measure}\label{EM}

\subsection{Some basic facts about equilibrium measure}

Let $D$ be a compact subset of the complex plane $ \mathbb{C}$, and $W(z)$
a potential function (a continuous function on $D$).
Denote by $\mathcal{M}(D)$ the set of the Borel
probability measures $\sigma$ on $D$ without point masses, and def\/ine the energy of
$\sigma$ to be
\begin{gather*}
I_\sigma=\int_D W(z)\text{d}\sigma(z)+\int_D\int_D
\log|z-w|^{-1}\text{d}\sigma(z)\text{d}\sigma(w).
\end{gather*}
An equilibrium measure for $W$ on $D$ is a
measure $\sigma\in \mathcal{M}(D)$ such that
\begin{gather*}
 I_\sigma=\inf_{\mu\in {\mathcal M}(D)}I_\mu.
 \end{gather*}

 \begin{theorem}\label{EU-equ-measure}
The equilibrium measure $\sigma$ exists and is unique.
It satisfies the equation
\begin{gather}\label{vary-form}
W(z)-2\int_D \log|z-w|\text{\rm d}\sigma(w)=C,
\end{gather}
where $C$ is a constant, almost everywhere with respect to
$\sigma$.
\end{theorem}

The proof of this theorem can be found in \cite{EF}.

Note that equation (\ref{vary-form}) does not have to hold
outside the support of $\sigma$.

Note also that if $\sigma$ is
absolutely continuous with respect to the Lebesgue measure
near a~point~$z_0$ in the interior of $D$, and
$\text{d}\sigma=g(z)\text{d}^2z$, where $g$ is continuous near $z_0$ and
$g(z_0)>0$, then $\Delta W=4\pi g$ near $z_0$. This clearly
cannot happen at points where $\Delta W\le 0$.
In particular, if $\Delta W\le 0$ everywhere,
then $\text{d}\sigma$ tends to be concentrated on the boundary of $D$.

\subsection{Asymptotic eigenvalue distribution in the normal matrix model}

In Section \ref{NMM}, we def\/ined
a measure
\begin{gather*}
P_{N}(M)\text{d} M=J_N(z_1,\ldots,z_N)\text{d}^2z_1\cdots\text{d}^2z_N\text{d} U.
\end{gather*}
by formula (\ref{NM-PF}).
We are interested in the behavior of this measure
when $N\to \infty$. Let
$\delta_{z}=\frac{1}{N}\sum\limits_{j=1}^{N}\delta_{z_j}$ be
the measure on $D$ corresponding to the points $z_j$. Then
\begin{gather*}
-\log (Z_NJ_N(z_1,\ldots,z_N))=N^2
\left(\int W(\xi)\text{d} \delta_{z}(\xi)+
\iint_{\xi\neq \zeta}\log|\xi-\zeta|^{-1}\text{d} \delta_{z}(\zeta)
\text{d} \delta_{z}(\xi)\right).
\end{gather*}
This shows that the leading contribution to the integral with
respect to the measure $P_N(M)\text{d} M$ comes from conf\/igurations of
eigenvalues $z_1,\ldots,z_N$ for which the expression in parentheses
in the last equation is minimized. This means that in the limit
$N\to \infty$, we should expect the measures $\delta_z$ for optimal conf\/igurations
to converge to the equilibrium measure with potential
function $W$.

This indeed turns out to be the case, as shown by the following
theorem, proved in \cite{EF}.

\begin{theorem}\label{limit}
Let the $k$-point correlation function be
\begin{gather*}
R_{N}^{(k)}((z_{i})_{i=1}^{k})=\int_{D^{N-k}}J_N(z_1,\ldots,z_N)
\prod_{i=k+1}^{N}\text{\rm d}^{2}z_{i}.
\end{gather*}
Then the measure
\begin{gather*}
R_{N}^{(k)}((z_{i})_{i=1}^{k})\prod_{i=1}^{k}\text{\rm d}^{2}z_{i}
\end{gather*}
on $D^{k}$ converges weakly to $\text{\rm d} \sigma^{\otimes k}$,
where $\text{\rm d} \sigma$ is the equilibrium measure on $D$, corresponding
to the potential function $W$.
\end{theorem}

In particular, if $k=1$, it means that the
eigenvalue distribution tends
to the equilibrium measure in $D$ as $N\to \infty$.

\subsection{Asymptotic eigenvalue distribution\\
in the generalized normal matrix model}

As we have seen above, the eigenvalue distribution in the
generalized normal matrix model is
\begin{gather*}
P_{N}(M)\text{d} M=J_{N,m}(z_1,\ldots,z_N)\text{d}^2z_1\cdots\text{d}^2z_N\text{d} U,
\end{gather*}
where
\begin{gather*}
-\log(2^{N}Z_NJ_{N,m})
=
N^2\left(\int W(\xi)\text{d} \delta_{z}(\xi)+
\iint_{\xi\neq \zeta}\log|\xi-\zeta|^{-1}\text{d} \delta_{z}(\zeta)
\text{d} \delta_{z}(\xi)\right)\\
\phantom{-\log(2^{N}Z_NJ_{N,m})=}{}
-\frac{N}{2}\int\log Q'(Q^{-1}(|\xi|^2))\text{d}\delta_z(\xi).
\end{gather*}
In the limit $N\to \infty$ the second term becomes unimportant
compared to the f\/irst one, which implies that Theorem
\ref{limit} is valid for the generalized normal matrix model.
Thus in the limit $N\to \infty$, the usual and the generalized
normal matrix models (with the same potential) are equivalent.

\section{Reconstruction of the boundary of the domain}\label{RBD}

In previous sections, we showed that in the normal matrix model
and the generalized normal matrix model, when $N\to \infty$, the
eigenvalue distribution converges to an equilibrium measure on
$D$ corresponding to some potential function $W$.
In this section, we will try to f\/ind this measure explicitly in
some special cases.

More specif\/ically, we will consider the case when $\Delta W>0$.
In this case, if the region $D$ is suf\/f\/iciently large,
it turns out that the equilibrium measure is often absolutely continuous with
respect to Lebesgue measure, and equals
$\text{d}\sigma=(4\pi)^{-1}\chi_E\Delta W\text{d}^2 z$, where $E$ is a region
contained in $D$ (the region of eigenvalues),
and $\chi_E$ is the characteristic function of
$E$. More precisely, it follows from Proposition 3.4
in \cite{EF} that if there exists a region $E\subset D$ such that
$\text{d}\sigma$ satisf\/ies equation (\ref{vary-form}) in $E$,
and the left hand side of this equation is $\ge C$
on $D\setminus E$, then $\text{d}\sigma$ is the equilibrium measure
on $D$ for the potential function $W$.
Moreover, note that if $E$ works for some $D$
then it works for any smaller $D'$ such that $E\subset D'\subset
D$. So, in a sense, $E$ is independent of $D$.
(Here we refer the reader to \cite{HM}, section 4, where
there is a much more detailed and precise treatment of
equilibrium measures, without the assumption $\Delta\Phi>0$).

Thus let us assume that $E$ exists, and consider
the problem of f\/inding it explicitly given the potential $W$.

\subsection{The reconstruction problem}

We will consider the case when $D=D(R)$ is the disk of radius $R$
centered at the origin, and
\begin{gather*}
W(z)=\Phi(z\bar z)-P(z)-\overline{P(z)},
\end{gather*}
where $\Phi$ is a function of one variable continuous on
$[0,\infty)$ and twice continuously dif\/ferentiable
on $(0,\infty)$, and $P$ a complex
polynomial. We assume that $(s\Phi'(s))'$ is positive, integrable
near zero, and satisf\/ies the boundary condition $\lim_{s\to 0}s\Phi'(s)=0$.
Computing the Laplacian of $W$,
we get (taking into account that $\Delta=4\partial \bar
\partial$):
\begin{gather*}
g(s):=(4\pi)^{-1}\Delta W=\pi^{-1}(\Phi'(s)+s\Phi''(s)),
\end{gather*}
where $s=z\bar z$. Def\/ine the measure
$\text{d}\sigma=g\text{d}^2z$.

Suppose that the region $E$ exists, and contains the origin.
In this case, dif\/ferentiating equation~(\ref{vary-form}) with
respect to $z$, we have inside $E$:
\begin{gather}\label{withp}
\bar z\Phi'(z\bar z)-P'(z)=\int_E \frac{g(w\bar w)}{z-w}\text{d}^2w.
\end{gather}
On the other hand, inside the disk $D$, the function
\begin{gather*}
W_0(z):=2\int_D g(w\bar w)\log|z-w|\text{d}^2w
\end{gather*}
satisf\/ies the equation $\Delta W_0=4\pi g$,
and is rotationally invariant, so
\begin{gather*}
W_0(z)=\Phi(z\bar z)+C',
\end{gather*}
where $C'$ is a constant. Hence, dif\/ferentiating, we get,
inside $D$:
\begin{gather}\label{nop}
\bar z\Phi'(z\bar z)=\int_D \frac{g(w\bar w)}{z-w}\text{d}^2w.
\end{gather}

Thus, subtracting (\ref{withp}) from (\ref{nop}), we obtain
inside $E$:
\begin{gather}\label{peq}
P'(z)=\int_{D\setminus E} \frac{g(w\bar w)}{z-w}\text{d}^2w.
\end{gather}

Let $I(s)=\pi\int_0^s g(t)\text{d} t=s\Phi'(s)$. Then
$\bar\partial I(z\bar z)=\pi z g(z\bar z)\text{d}\bar z$.
Thus, using Green's formula, we get from (\ref{peq}):
\begin{gather*}
P'(z)=\frac{1}{2\pi i}
\int_{\partial D-\partial E} \frac{I(w\bar w)}{w(z-w)}\text{d} w,
\end{gather*}
where the boundaries are oriented counterclockwise.
The integral over the boundary of $D$ is zero by Cauchy's formula,
so we are left with the equation
\begin{gather*}
P'(z)=\frac{1}{2\pi i}\int_{\partial E} \frac{I(w\bar w)}{w(w-z)}\text{d} w.
\end{gather*}

This equation appeared f\/irst in the theory of Hele-Shaw
f\/lows on curved surfaces in \cite{EE}, and it can be solved
explicitly by the method of singular points developed in the same
paper. Let us recall this method.

\subsection{The singular point method}
Def\/ine the Cauchy transform $h_E$ of $E$ with respect to the measure
$\text{d}\sigma$ by
\begin{gather*}
h_E(z)=\int_{D\setminus E}\frac{\text{d}\sigma(w)}{z-w},\qquad  z\in E.
\end{gather*}
This is a holomorphic function of $z$ which (as we have just
seen) is independent of the radius $R$ of $D$. As we have seen,
it is also given by the contour integral
\begin{gather}\label{contint}
h_E(z)=\frac{1}{2\pi i}\int_{\partial E} \frac{I(w\bar
w)}{w(w-z)}\text{d} w,
\end{gather}
and in our case we have $h_E(z)=P'(z)$.

Let $f: D(1)\to \overline{\Bbb C}\setminus E$
be a conformal map, such that $f(0)=\infty$, and $(1/f)'(0)=a\in
\Bbb R^+$ (such a map is unique).

\begin{lemma}
The function
\begin{gather*}
\phi(\zeta)=I(f(\zeta)\overline{f(\zeta)})-f(\zeta)h_E(f(\zeta))
\end{gather*}
continues analytically from the unit circle
to a holomorphic function outside the unit disk.
\end{lemma}

\begin{proof}
By the Cauchy formula, we have
\begin{eqnarray*}
h_E(z)=\frac{1}{2\pi i}\int_{\partial E}\frac{h_E(w)}{w-z}\text{d} w,\qquad
\text{for any} \ \ z\in E.
\end{eqnarray*}

So by formula (\ref{contint}), we have
\begin{gather*}
\frac{1}{2\pi i}\int_{\partial
E}\frac{I(w\bar{w})/w-h_E(w)}{w-z}\text{d} w=0,\qquad
\text{for any} \ \ z\in E.
\end{gather*}
It follows that the function $I(z\bar{z})/z-h_E(z)$,
def\/ined along $\partial E$,
can be analytically continued to a holomorphic function outside
$E$, which vanishes at inf\/inity. This implies the lemma.
\end{proof}

Similarly to \cite{EE}, this lemma implies
the following theorem.

\begin{theorem}\label{reconstruction}
The function $h_E$ is rational if and only if the function
$\theta(\zeta)=I(f(\zeta)\overline{f(1/\bar{\zeta})})$ is.
Moreover, the number of poles of $\theta$ is twice of the
number of poles of $zh_E(z)$. More specifically,
if $\zeta_0$ and $1/\bar\zeta_0$ are poles of $\theta$ of order
$m$, then $z_0=f(\zeta_0)$ is a pole of order $m$ for $h_E(z)$,
and vice versa.
\end{theorem}

Thus, if $h$ is a rational function, then $\theta$ can be determined from
$h$ up to f\/initely many parameters.

After this, $f$ can be reconstructed from $\theta$ using the
Cauchy formula. For this, note that the function $I$ is
invertible, since $I'=g>0$. Also, $\theta$ takes nonnegative real
values on the unit circle. Thus, we have
\begin{gather*}
f(\zeta)\overline{f(1/\bar \zeta)}=I^{-1}(\theta(\zeta)).
\end{gather*}
Taking the logarithm of both sides, we obtain
\begin{gather*}
\log (\zeta f(\zeta))+\log \big(\zeta^{-1}\overline{f(1/\bar
\zeta)}\big)=\log I^{-1}(\theta(\zeta)).
\end{gather*}
Thus we have
\begin{gather*}
f(\zeta)= a\zeta^{-1}\exp\left(\frac{1}{2\pi i}
\int_{|u|=1}\frac{\log I^{-1}(\theta(u))}{u-\zeta}\text{d} u\right),\\
a= \exp\left(-\frac{1}{4\pi i}\int_{|u|=1}\frac{\log
I^{-1}(\theta(u))}{u}\text{d} u\right).
\end{gather*}

The unknown parameters of $\theta$ can now be determined
from the cancellation of poles in Theorem \ref{reconstruction},
similarly to the procedure described in \cite{VE}.
We note that the knowledge of the function $h_E$ is not suf\/f\/icient
to determine $E$ (for example if $E$ is a disk of any radius
centered at the origin then $h_E=0$). To determine the parameters
completely, we must also use the information on the area of $E$:
\begin{gather*}
\int_E \text{d}\sigma=-\frac{1}{2\pi i}
\int_{|u|=1}\frac{\theta(u)}{f(u)}f'(u)\text{d} u=1.
\end{gather*}

\subsection{The polynomial case}
In particular, in our case,
\begin{gather*}
h_E(z)=P'(z)=a_1+a_2z+\dots+a_d z^{d-1},
\end{gather*}
which implies that
$\theta(\zeta)=\sum\limits_{j=-d}^{d} b_j\zeta^j$, and $\bar b_j=b_{-j}$.

So we get
\begin{gather}
f(\zeta)=a\zeta^{-1}\exp\left(\frac{1}{2\pi i}
\int_{|u|=1}\frac{\log I^{-1}\left(\sum\limits_{j=-d}^{d} b_ju^j\right)}
{u-\zeta}\text{d} u\right),\nonumber\\
a=\exp\left(-\frac{1}{4\pi i}\int_{|u|=1}\frac{\log
I^{-1}\left(\sum\limits_{j=-d}^{d} b_ju^j\right)}{u}\text{d} u\right).\label{conformal-map1}
\end{gather}

Finally, note that if the coef\/f\/icients of the polynomial $P$ are
small enough, then all our assumptions are satisf\/ied:
the region $E$ exists (in fact, it is close to a disk),
and contains the origin. Also, in this case the left hand side
of equation (\ref{vary-form}) is $\ge C$, which implies
that the equilibrium measure in this case (and hence, the
asymptotic eigenvalue distribution) is the measure $\text{d} \sigma$ in
the region $E$.

\begin{example}
Consider Example \ref{mmdag}:
the generalized normal matrix model with the
density $\exp(-\beta\textbf{tr}(MM^\dagger-P(M)-P(M)^\dagger)$.
As we showed, in this case
$\Phi(s)=m\beta Q^{-1}(s)$.
So a~short computation shows that
\begin{gather*}
\pi m^{-1}\beta^{-1} g(Q(x))=\frac{Q(x)^2}{Q'(x)^3}\sum_i
\frac{1}{(x+\alpha_i)^2}.
\end{gather*}
This implies that $g>0$, i.e.\ our analysis
applies in this case.
\end{example}

\subsection{Some explicit solutions}

Consider the case $\Phi(s)=Cs^b$, $C,b>0$.
For example, in the generalized normal matrix
model with $\alpha_i=0$ and potential term
as in Example \ref{mmdag}, one has $\Phi(s)=ms^{1/m}$,
which is a special case of the above.

We have $g(s)=\pi^{-1}Cb^2s^{b-1}$, so our analysis
applies (note that if $b<1$ then $g$ is singular at zero, but the
singularity is integrable and thus nothing really changes in
our considerations), and $I(s)=Cbs^b$. Thus the integral in
(\ref{conformal-map1}) can be computed explicitly (by factoring
$\theta$), and the formula for the conformal map
$f$ simplif\/ies as follows:
\begin{gather*}
f(\zeta)=(a\zeta)^{-1}\prod_{j=1}^d (1-\zeta \zeta_j^{-1})^{1/b}.
\end{gather*}

The parameters $a>0$ and $\zeta_j$ are determined from
the singularity conditions and the area condition.

Consider for simplicity the example $d=1$.
In this case we have
\begin{gather*}
h_E(z)=K,
\end{gather*}
and we can assume without loss of generality that $K\in \Bbb R$.
Then
\begin{gather*}
f(\zeta)=(a\zeta)^{-1}(1+\beta\zeta)^{1/b},\qquad \beta\in \Bbb R,
\end{gather*}
and
\begin{gather*}
\theta(\zeta)=Cba^{-2b}(1+\beta\zeta)(1+\beta\zeta^{-1}).
\end{gather*}
The residue of $\theta$ at zero is thus $Cb\beta a^{-2b}$.
Thus the singularity condition says
\begin{gather*}
Cb\beta a^{1-2b}=K.
\end{gather*}
The area condition is
\begin{gather*}
1=Cba^{-2b}(1+\beta^2(1-b^{-1})).
\end{gather*}
Thus we f\/ind
$
\beta=KC^{-1}b^{-1}a^{2b-1}$,
and the equation for $a$ has the form
\begin{gather*}
Cba^{-2b} + C^{-1}b^{-1}K^2a^{2b-2}(1-b^{-1})=1.
\end{gather*}

\begin{remark} This example shows that to explicitly
solve the generalized (as opposed to the usual)
normal matrix model in the $N\to \infty$ limit
with the quadratic (Gaussian) potential, one really needs
the technique explained in Section 5 of this paper,
and the techniques of \cite{EF} are not suf\/f\/icient.
\end{remark}

\subsection*{Acknowledgements}
P.E. is grateful to G.~Felder and
P.~Wiegmann for useful discussions.
The work of P.E. was  partially supported by the NSF grant
 DMS-0504847.

\pdfbookmark[1]{References}{ref}

\LastPageEnding

\begin{thebibliography}{99}

\footnotesize\itemsep=0pt

\bibitem{CZ}
Chau L.-L., Zaboronsky O.,
On the structure of correlation functions in the normal matrix model,
{\it Comm. Math. Phys.} {\bf 196} (1998), 203--247, \href{http://arxiv.org/abs/hep-th/9711091}{hep-th/9711091}.

\bibitem{CBH} Crawley-Boevey W., Holland M.P.,
Noncommutative deformations of Kleinian singularities,
{\it Duke Math. J.} {\bf 92} (1998), 605--635.

\bibitem{EF}
Elbau P., Felder G.,
Density of eigenvalues of random normal matrices,
{\it Comm. Math. Phys.} {\bf 259} (2005), 433--450, \href{http://arxiv.org/abs/math.QA/0406604}{math.QA/0406604}.

\bibitem{EE}
Entov V.M., Etingof P.I.,
Viscous f\/lows with time-dependent free boundaries in a non-planar Hele-Shaw cell,
{\it Euro. J. Appl. Math.} {\bf 8} (1997), 23--35.

\bibitem{HM}
Hedenmalm H., Makarov N.,
Quantum Hele-Shaw f\/low,
\href{http://arxiv.org/abs/math.PR/0411437}{math.PR/0411437}.

\bibitem{KKMWZ}
Kostov I.K., Krichever I., Mineev-Weinstein M., Wiegmann P.B., Zabrodin A.,
The $\tau$-function for analytic curves, in Random Matrix Models and Their Applications,
{\it Math. Sci. Res. Inst. Publ.}, Vol.~40, Cambridge Univ. Press, Cambridge, 2001, 285--299.

\bibitem{KMZ}
Krichever I., Marshakov A., Zabrodin A.,
Integrable structure of the Dirichlet boundary problem in multiply-connected domains, {\it Comm. Math. Phys.} {\bf 259} (2005), 1--44, \href{http://arxiv.org/abs/hep-th/0309010}{hep-th/0309010}.

\bibitem{MWZ}
Marshakov A., Wiegmann P.B., Zabrodin A.,
Integrable structure of the Dirichlet boundary problem in two dimensions,
{\it Comm. Math. Phys.} {\bf 227} (2002), 131--153, \href{http://arxiv.org/abs/hep-th/0109048}{hep-th/0109048}.

\bibitem{O}
Oas G.,
Universal cubic eigenvalue repulsion for random normal matrices,
{\it Phys. Rev. E} {\bf 55} (1997), 205--211, \href{http://arxiv.org/abs/cond-mat/9610073}{cond-mat/9610073}.

\bibitem{VE}
Varchenko A.N., Etingof P.I.,
Why the boundary of a round drop becomes a curve of order four,
AMS, Providence, 1992.

\bibitem{WZ1}
Wiegmann P.B., Zabrodin A.,
Conformal maps and integrable hierarchies,
{\it Comm. Math. Phys.} {\bf 213} (2000), 523--538, \href{http://arxiv.org/abs/hep-th/9909147}{hep-th/9909147}.

\bibitem{WZ}
Wiegmann P.B., Zabrodin A.,
Large scale correlations in normal non-Hermitian matrix ensembles,
{\it J.~Phys.~A: Math. Gen.} {\bf 36} (2003), 3411--3424,
\href{http://arxiv.org/abs/hep-th/0210159}{hep-th/0210159}.

\end{thebibliography}
\end{document}